\documentclass[11pt]{article}
\usepackage{amsfonts}
\usepackage{amssymb}
\usepackage{epsfig}
\usepackage{psfrag}
\usepackage{url}  

\def\@begintheorem#1#2{\par\bgroup{\sc #1\ #2. }\it\ignorespaces}
\def\@opargbegintheorem#1#2#3{\par\bgroup{\sc #1\ #2\ (#3). } \it\ignorespaces}
\def\@endtheorem{\egroup}
\newtheorem{theorem}{Theorem}[section]
\newtheorem{corollary}[theorem]{Corollary}
\newtheorem{lemma}[theorem]{Lemma}
\newtheorem{proposition}[theorem]{Proposition}
\newtheorem{problem}[theorem]{Problem}
\newtheorem{example}[theorem]{Example}
\newtheorem{remark}[theorem]{Remark}
\newtheorem{algorithm}[theorem]{Algorithm}
\newtheorem{definition}[theorem]{Definition}
\newtheorem{conjecture}[theorem]{Conjecture}
\newcommand{\bt}[1]{\begin{theorem}\label{#1}}
\newcommand{\bc}[1]{\begin{corollary}\label{#1}}
\newcommand{\bl}[1]{\begin{lemma}\label{#1}}
\newcommand{\bp}[1]{\begin{proposition}\label{#1}}
\newcommand{\bpro}[1]{\begin{problem}\label{#1}}
\newcommand{\be}[1]{\begin{example}\rm\label{#1}}
\newcommand{\ba}[1]{\begin{algorithm}\rm\label{#1}}
\newcommand{\bd}[1]{\begin{definition}\rm\label{#1}}
\newcommand{\bpr}{\par \noindent{\it Proof}. \ignorespaces}
\newcommand{\bprof}[1]{\par\noindent{\it Proof of #1}. \ignorespaces}
\newcommand{\et}{\end{theorem}}
\newcommand{\ec}{\end{corollary}}
\newcommand{\el}{\end{lemma}}
\newcommand{\ep}{\end{proposition}}
\newcommand{\epro}{\end{problem}}
\newcommand{\ee}{\end{example}}
\newcommand{\ea}{\end{algorithm}}
\newcommand{\ed}{\end{definition}}
\newcommand{\epr}{{\ \vbox{\hrule\hbox{
\vrule height1.3ex\hskip0.8ex\vrule}\hrule
}}\\\par}
\newcommand{\reals}{\mathbb R}
\newcommand{\rationals}{\mathbb Q}
\newcommand{\cone}{\hbox{cone}}

\def\R{\mathbb{R}}

\begin{document}

\title{\bf Graphs of Transportation Polytopes}

\author{
Jes\'us A. De Loera
\thanks{Research supported in part by NSF grant DMS-0608785.}
\and
Edward D. Kim
\thanks{Research supported in part by VIGRE NSF grant DMS-013534.}
\and
Shmuel Onn
\thanks{Supported in part by a grant from ISF - the
Israel Science Foundation, by the Technion President Fund,
and by the Fund for the Promotion of Research at the Technion.}
\and
Francisco Santos
\thanks{Supported in part by the Spanish Ministry of Science through grant 
MTM2008-04699-C03-02, and a mobility grant}
}

\date{}
\maketitle

\begin{abstract}
  This paper discusses properties of the graphs of $2$-way and $3$-way
  transportation polytopes, in particular, their possible numbers of
  vertices and their diameters. Our main results include a quadratic
  bound on the diameter of axial $3$-way transportation polytopes and
  a catalogue of non-degenerate transportation polytopes of small
  sizes. The catalogue disproves five conjectures about these
  polyhedra stated in the monograph by Yemelichev et al. (1984). It
  also allowed us to discover some new results.  For example, we prove
  that the number of vertices of an $m\times n$ transportation
  polytope is a multiple of the greatest common divisor of $m$ and
  $n$.
\end{abstract}

\section{Introduction}

This paper takes a new look at the graphs of transportation polytopes.
Transportation polytopes are well-known objects in operations research
and mathematical programming (see e.g., \cite{BB, BR,barneshoffman,
  CDMS, KW, OR, OS, Vla, YKK} and references therein).  Statisticians
have also interest in them of their own (see e.g., \cite{Cox1, Cox2,
  DG, DFKPR, IJ, MP} and references therein).

During the 1970's and 1980's the study of the classical $2$-way
transportation polytopes, i.e., those polytopes of $m \times n$ tables
satisfying row-sum and column-sum conditions, was very active. The
state of the art of that research was carefully summarized in the
comprehensive book by Yemelichev, Kovalev, and Kratsov \cite{YKK} and
Vlach's survey \cite{Vla}.  On the other hand, 3-way transportation
polytopes, whose feasible points are $l \times m \times n$ arrays of
non-negative numbers satisfying certain sum conditions, are less
understood.  We define the polyhedra whose points are $l \times m
\times n$ tables satisfying certain sum conditions. They come in two
main varieties:
\begin{enumerate}
\item First, consider the \emph{axial} 3-way 
transportation polytope: Let $x=(x_1,\dots,x_l)$,
$y=(y_1,\dots,y_m)$, and $z=(z_1,\dots,z_n)$ be three rational vectors
of lengths $l$, $m$ and $n$, respectively, with non-negative entries.
Let $T_{x,y,z}$ be the polytope defined by \emph{1-marginals}; that
is, the following $l+m+n$ equations in the $l\times m\times n$ real
variables $a_{ijk}$ ($i=1,\dots,l; j=1,\dots, m; k=1,\dots, n$):
\[ 0\le a_{i,j,k}, \forall i,j,k\qquad \sum_{j,k} a_{i,j,k} = x_{i}, \forall
i\qquad \sum_{i,k} a_{i,j,k} = y_{j}, \forall j\qquad \sum_{i,j}
a_{i,j,k} = z_{k}, \forall k.  \label{eq:1marginals} \]

Observe that a necessary and sufficient condition for $T_{x,y,z}$ to 
be non-empty is that
\[
\sum_i x_i = \sum_j y_j = \sum_k z_k,
\label{eq:axialsizigies}
\]
and that, consequently, $T_{x,y,z}$ is defined by only $l+m+n-2$
independent equations.
\item Similarly, \emph{planar} $3$-way transportation polytopes can be
  defined by specifying three matrices $U\in M_{m,n}(\rationals) ,\
  V\in M_{l,n}(\rationals),\ W\in W_{l,m}(\rationals)$ for the
  line-sums resulting from fixing two of the indices of entries and
  adding over the remaining index. That is, we have the following
  $lm+ln+mn$ equations (\emph{2-marginals}) in the same $lmn$
  variables $a_{i,j,k}$:
\[
0\le a_{i,j,k}, \forall i,j,k\qquad
\sum_{i} a_{i,j,k} = U_{j,k}, \forall j,k\qquad 
\sum_{j} a_{i,j,k} = V_{i,k}, \forall i,k\qquad 
\sum_{k} a_{i,j,k} = W_{i,j}, \forall i,j.
\label{2marginals}
\]
One can see that in fact only $lm + ln + mn - l - m - n + 1$ of the
defining equations are linearly independent for feasible systems.
\end{enumerate}

Observe that the axial 3-way transportation polytopes generalize the
classical transportation polytope of size $m\times n$, which coincides
with $T_{x,y,z}$ for $l=1$ and $x=\sum y_j=\sum z_k$. A less trivial
rewriting of the classical $2\times n$ transportation polytope as a
$2\times 2\times n$ 3-way planar transportation polytope is given in
Theorem~\ref{thm:22n} below.

Recall that the \emph{$1$-skeleton} or \emph{graph} of a convex
polytope $P$ is the set of all $0$-dimensional and $1$-dimensional
faces (vertices and edges) of $P$, with their natural incidence
relation. The main focus of this paper is to investigate the number of
vertices and the diameters of the graphs of classical (that is,
$2$-way) and $3$-way transportation polytopes.

Some of the statements below require our transportation polytopes to
be \emph{non-degenerate}.  By this we mean that the polytope is simple
(i.e., every vertex is adjacent to dimension many other vertices) and
it is of maximal possible dimension (that is, dimension $lmn-l-m-n+2$
for $l \times m \times n$ axial transportation polytopes, and
dimension $(l-1)(m-1)(n-1)$ for $l \times m \times n$ planar
transportation polytopes). Graphs of non-degenerate transportation
polytopes are of particular interest because they have the largest
possible number of vertices and largest possible diameter among the
graphs of all transportation polytopes of given type and parameters.
Indeed, if a transportation polytope $P$ is degenerate, by carefully
perturbing the marginals that define it we can get a non-degenerate
one $P'$.  The perturbed marginals are obtained by taking a feasible
point $X$ in $P$, perturbing the entries in the table and using the
recomputed sums as the new marginals for $P'$.  The graph of $P$ can
be obtained from that of $P'$ by contracting certain edges, which
cannot increase either the diameter nor the number of vertices.

Our main result is a bound on the diameter of axial 3-way
transportation polytopes:
\begin{theorem}
\label{theorem:main}
The graph of the $l \times m \times n$ axial transportation polytope
$T_{x,y,z}$ has diameter at most $2(l+m+n-3)^2$.
\end{theorem}

A similar result for the graph of a classical transportation polytope
was given by Brightwell et al. (see \cite{BvHS}), who proved an upper
bound of $8(m+n-2)$ for the diameter.  More recently, Hurkens (see 
\cite{hurkens}) has obtained a bound of $4(m+n-1)$, a factor of four
away from the predicted value of the Hirsch conjecture.

Using computational tools, we also give a complete catalogue of
\emph{non-degenerate} $2$-way and $3$-way transportation polytopes
(both axial and planar) of small sizes. This allowed us to explore
properties of transportation polytopes (e.g., their diameters and how
close they were to the Hirsch conjecture bound).  The summary of the
catalogue is: 
\begin{theorem}\label{genericcase}
\begin{itemize}

\item The only possible numbers of vertices of non-degenerate $2
  \times 3$, $2 \times 4$, $2 \times 5$, $3 \times 3$, and $3 \times
  4$ classical transportation polytopes are those given in
  Table~\ref{classical_vertices}.

\item The only possible numbers of vertices of non-degenerate $2
  \times 2 \times 2$, $2 \times 2 \times 3$, $2 \times 2 \times 4$, $2
  \times 2 \times 5$, and $2 \times 3 \times 3$ planar transportation
  polytopes are those given in Table \ref{planar_vertices}.

  Every non-degenerate $2 \times 3 \times 4$ planar transportation
  polytope has between $7$ and $480$ vertices.

\item The only possible numbers of vertices of non-degenerate $2
  \times 2 \times 2$ and $2 \times 2 \times 3$ axial transportation
  polytopes are those given in Table \ref{axial_vertices}.

  Every non-degenerate $2 \times 2 \times 4$ axial transportation
  polytope has between $32$ and $504$ vertices.  Every non-degenerate
  $2 \times 3 \times 3$ axial transportation polytopes has between
  $81$ and $1056$ vertices.  The number of vertices of non-degenerate
  $3 \times 3 \times 3$ axial transportation polytopes is at least
  $729$.
\end{itemize}
\end{theorem}

\begin{table}[hbtp]
 \font\ninerm=cmr8
 \centerline{\ninerm
 \begin{tabular}{|c|c|c|}
 \hline
\textrm{\qquad Size\qquad} & \textrm{Dimension} & \textrm{Possible numbers of vertices} \cr \hline         
$2 \times 3$ & $2$ & 3 4 5 6 \cr \hline
$2 \times 4$ & $3$ & 4 6 8 10 12 \cr \hline
$2 \times 5$ & $4$ & 5 8 11 12 14 15 16 17 18 19 20 21 22 23 24 25 26 27 28 29 30 \cr \hline
$3 \times 3$ & $4$ & 9 12 15 18 \cr \hline
$3 \times 4$ & $6$ & 16 21 24 26 27 29 31 32 34 36 37 39 40 41 \cr
             &     & 42 44 45 46 48 49 50 52 53 54 56 57 58 60 61 \cr
             &     & 62 63 64 66 67 68 70 71 72 74 75 76 78 80 84 90 96 \cr \hline
\end{tabular}
}
\caption{Numbers of vertices possible in non-degenerate classical transportation polytopes} \label{classical_vertices}
\end{table}

\begin{table}[hbtp]
 \font\ninerm=cmr8
 \centerline{\ninerm
 \begin{tabular}{|c|c|c|}
 \hline
\textrm{\qquad Size\qquad} & \textrm{Dimension} & \textrm{Possible numbers of vertices} \cr \hline         
$2 \times 2 \times 2$ & $1$ & 2 \\ \hline
$2 \times 2 \times 3$ & $2$ & 3 4 5 6 \\ \hline
$2 \times 2 \times 4$ & $3$ & 4 6 8 10 12 \\ \hline
$2 \times 2 \times 5$ & $4$ & 5 8 11 12 14 15 16 17 18 19 20 21 22 23 24 25 26 27 28 29 30 \\ \hline
$2 \times 3 \times 3$ & $4$ & 5 8 9 11 12 13 14 15 16 17 18 19 20 21 22 23 \cr
                      &     & 24 25 26 27 28 29 30 31 32 33 34 35 36 37 38 39 40 \cr
                      &     & 41 42 43 44 45 46 47 48 49 50 51 52 53 54 55 56 57 58 59 \\ \hline
\end{tabular}
}
\caption{Numbers of vertices possible in non-degenerate planar transportation polytopes} \label{planar_vertices}
\end{table}

\begin{table}[hbtp]
 \font\ninerm=cmr8
 \centerline{\ninerm
 \begin{tabular}{|c|c|c|}
 \hline
\textrm{\qquad Size\qquad} & \textrm{Dimension} & \textrm{Possible numbers of vertices} \cr \hline
$2 \times 2 \times 2$ & $4$ & 8 11 14 \\ \hline
$2 \times 2 \times 3$ & $7$ & 18 24 30 32 36 38 40 42 44 46 48 50 52 54 56 58 60 62 64 66 \cr
                      &     & 68 70 72 74 76 78 80 84 86 96 108 \\ \hline
\end{tabular}
}
\caption{Numbers of vertices possible in non-degenerate axial transportation polytopes} \label{axial_vertices}
\end{table}

The catalogue was obtained via the exhaustive and systematic computer
enumeration of all combinatorial types of non-degenerate
transportation polytopes. The theoretical foundations of it are the
notions of parametric linear programming, chamber complex, Gale
diagrams and secondary polytopes (see \cite{DRS}). We present them in
Section~\ref{foundations}.  The full catalogue of transportation
polytopes (including other families, such as $3 \times 5$, $4 \times
4$, $4 \times 5$, etc.) is available in a searchable web database at:
\url{http://www.math.ucdavis.edu/~ekim/transportation_polytope_database/}.

Based on the data we collected, we discovered and proved (in
Section~\ref{sec:gcd}) the following results:
\begin{theorem}
\label{thm:gcd}
The number of vertices of a non-degenerate $m \times n$ classical
transportation polytope is divisible by $\hbox{GCD}(m,n)$.
\end{theorem}

\begin{theorem}
\label{thm:22n}
The $2 \times 2 \times n$ planar transportation polytopes are in $1-1$
correspondence with the $2 \times n$ classical transportation
polytopes, with corresponding pairs being linearly isomorphic.
\end{theorem}

Note that Theorem \ref{thm:22n} is best possible in the sense that for
$m,n\geq 3$ there are many more types of planar $2\times m\times n$
transportation polytopes than types of $m\times n$ transportation
polytopes - see the rows of the $3\times 3$ and $2\times 3\times 3$
polytopes in Tables \ref{classical_vertices} and
\ref{planar_vertices}. Finally, we state two open conjectures at the
end of Section~\ref{foundations}.

We close the introduction explaining the relevance of our results to
current research. Bounding the diameter of graphs of polytopes has
received a lot of attention because of its connection to the
performance of the simplex method for linear programming and, most
especially, to try to understand the Hirsch conjecture (see
\cite{KO,KK} and references therein).

We recall that the Hirsch conjecture asserts that the diameter of the
graph of any polytope of dimension $d$ and with $f$ facets is bounded
above by $f-d$. Not only is this conjecture open, but even the weaker
statement asserting a polynomial upper bound (in $f$ and $d$) for the
diameters of graphs of polytopes is unknown (although a
quasi-polynomial bound appeared in \cite{KK}).

As we observed, \cite{BvHS} provided the first linear bound for the
diameter of the graphs of $2$-way transportation polytopes.
Theorem~\ref{theorem:main} provides a quadratic one for axial 3-way
transportation polytopes and, moreover, a sublinear one if we assume
that the three parameters $l$, $m$ and $n$ are approximately the same.
(Observe that the number of facets of a 3-way transportation polytope
is bounded above by the product $lmn$ of its size parameters).

Bounding the diameter of $3$-way transportation polytopes is
particularly interesting because of the following results recently
proved by two of the authors in \cite{DO2}:
\begin{enumerate}
\item Any rational convex polytope can be rewritten as a face $F$ of
  an axial $3$-way transportation polytope.  The sizes $l,m,n$,
  $1$-marginals $x,y,z$, and the entries $a_{i,j,k}$ that are
  prescribed to be zero in the face $F$ can be computed in polynomial
  time on the size of the input.

\item More dramatically, any convex rational polytope is
  isomorphically representable as a \emph{planar} 3-way transportation
  polytope. 
\end{enumerate}

That is to say, a version of Theorem~\ref{theorem:main} for the 3-way
planar case, or a version for the axial case that allows one to
prescribe some variables to be zero, would provide a polynomial upper
bound on the diameter of the graph of \emph{every} convex rational
polytope.

Another consequence of these results is that the method of   Section~\ref{foundations}
for enumerating all  combinatorial types of planar $3$-way  transportation polytopes, 
yields, in particular, an enumeration of all  types of
  rational convex polytopes.

Let us finally mention that our systematic listing of non-degenerate
transportation polytopes provides the solution to at least four open
problems and conjectures about transportation polytopes stated in the
monograph \cite{YKK}:

\begin{enumerate}

\item Klee and Witzgall in \cite{KW} prove that the largest possible
  number of vertices in classical transportation polytopes of size
  $m\times n$ is achieved by the {\em generalized Birkhoff polytope}
  (the transportation polytope with parameters $x_i=n\ \forall\ i$,
  $y_j=m\ \forall\ j$).  Problem 32 in page 400 of \cite{YKK}
  conjectured that the same holds in general.

  But in Example~\ref{exm:notBirkhoff} we provide an explicit
  counterexample of this for planar 3-way transportation polytopes.
  (In this case, the {\em generalized Birkhoff polytope} is the planar
  $3$-way transportation polytope whose 2-marginals are given by the
  $m \times n$ matrix $U(j,k)=l$, the $l \times n$ matrix $V(i,k)=m$,
  and the $l \times m$ matrix $W(i,j)=n$.)

\item Question 36 on page 396 of \cite{YKK} asked: \emph{Is it true
    that every integer of the form $(l-1)(m-1)(n-1)+t$ where $1 \leq t
    \leq ml+nl+mn-l-m-n$, and only these integers, can equal the
    number of facets of a non-degenerate planar $3$-way transportation
    polytope of order $l \times m \times n$, where $l,m,n \geq 2$?}

  For the case $l=m=2$ and $n=3$, the conjecture asks if every integer
  from $3$ to $11$, and only these integers, equal the number of
  facets of non-degenerate $2 \times 2 \times 3$ planar transportation
  polytopes. Since in this case the number of facets equals the number
  of vertices (because the polytopes are two dimensional)
  Table~\ref{planar_vertices} answers the question negatively: only
  facet-counts from 3 to 6 occur, while 7 through 11 are in fact
  missing.

\item Similarly, Conjecture 33 on page 400 of \cite{YKK} asked:
  \emph{Is it true that every integer from $1$ to $ml+nl+mn-l-m-n+1$,
    and only these numbers, are realized as the diameter of a planar
    $3$-way transportation polytope of order $l \times m \times n$?}

  The same case $l=m=2$, and $n=3$ shows that this is false. The
  transportation polytopes obtained are polygons with up to six sides,
  hence of diameter at most three, instead of 10.

\item Open problem 37 in page 396 of \cite{YKK} asks \emph{whether the
    numbers of vertices of $l \times m \times n$ non-degenerate planar
    transportation polytopes satisfy:}
\[
(l-1)(m-1)(n-1)+1 < f_0 < 2(l-1)(m-1)(n-1).
\]
We show the answer is no even in the case  $2\times 2\times 4$.

\end{enumerate}

In addition to the four solved problems above, Theorems
\ref{genericcase} and \ref{thm:gcd} are initial steps on the solution
of Problem 25 in page 399 of \cite{YKK}. It asks to find the complete
distribution of possible number of vertices for transportation
polytopes.

\section{Classifying Transportation Polytopes} \label{foundations}

Theorem \ref{genericcase} was obtained through an exhaustive
enumeration whose foundation is the theory of secondary polytopes and
parametric linear programming. In some cases when the full enumeration
was impossible we at least get lower and upper bounds for the number
of vertices that these polytopes can have. In this section we discuss
the necessary background to understand the construction of the
complete catalogue.

\subsection{Enumeration via Regular triangulations and Secondary polytopes}

We begin by recalling some basic facts about convex polytopes
presented, as all transportation polytopes are, in the form $P_c = \{x
: Bx=c, x \geq 0 \}$.  For the case of transportation polytopes, the
vector $c$ is the vector given by the demand/supply quantities. Fix a
matrix $B$ of full row rank.  Most of our results are obtained by
studying what happens to the combinatorics of $P_c$ as the vector $c$
changes while we fix the matrix $B$. This study, for general matrices,
is known as \emph{parametric linear programming} (see Chapter 1 of
\cite{DRS}).

A subset of $\R^n$ that is closed under addition and under
multiplication by positive scalars is a {\em cone}.  For any set $L$
of vectors in $\R^n$, the {\em cone generated by $L$}, denoted
$\cone(L)$, is the set of all vectors that can be expressed as
non-negative linear combinations of the members of $L$. Abusing
notation, for a matrix $B$, by $\cone(B)$ we mean the cone generated
the set of column vectors of $B$.

A maximal linearly independent subset $b$ of $B$ is a \emph{basis} of $B$.
Geometrically, each basis of the matrix $B$ spans a simple cone inside
$\cone(B)$.  Every basis $b$ of $B$ defines a \emph{basic solution} of
the system as the unique solution of the $m$ linearly independent
equations $bx_b=c$ and $x_j=0$ for $j$ not in $b$. A basic solution is
\emph{feasible} if in addition, $x\geq 0$. Geometrically, a basic
feasible solution corresponds to a simple cone that contains $c$. In
fact, one can see that $P_c$ is non-empty if and only if $c\in
\cone(B)$.

A fundamental fact in linear programming is that, for a given
right-hand-side vector $c$, all vertices of the polyhedron $P_c$ are
basic feasible solutions (see \cite{schrijver,YKK}). Moreover, if the
polyhedron $P_c$ is assumed to be non-degenerate then the basic
feasible solution must be strictly positive on the entries
corresponding to the basis $b$.  Geometrically, a basis $b\subseteq B$
produces a vertex of a non-degenerate $P_c$ if and only if $c$ lies in
the interior of the cone generated by $b$. In conclusion, the vertices
of a non-degenerate polyhedron $P_c$, from this family of parametric
polytopes, are in bijection with the bases that contain the
right-hand-side vector $c$ within their interior.

We now look at what happens when we let $c$ vary. If $P_c$ is
non-degenerate and the change in $c$ is small, the facets of $P_c$
move but the combinatorial type of $P_c$ does not change. Only when a
basic solution changes from being feasible to not feasible, or
viceversa, the combinatorics of $P_c$ (that is, the face lattice and,
in particular, the graph of $P_c$) can change.

Put differently: Denote by $\Sigma_B$ the set of all cones generated
by bases of $B$. Let $\partial \Sigma_B$ denote the union of the
boundaries of all elements of $\Sigma_B$. The connected components of
$\cone(B) \setminus \partial \Sigma_B$ are open convex cones called
the {\em chambers} of $B$.  Equivalently, we call the \emph{chamber
  associated to} a given right-hand side vector $c$ the intersection
of the interiors of simple cones that contain $c$ in their interior.
We remark that every feasible and sufficiently generic $c$ is in a
chamber (as opposed to lying on $\partial \Sigma_B$).  Two vectors
$c_1$ and $c_2$ in the same chamber determine non-degenerate polytopes
$P_{c_1}$ and $P_{c_2}$ that are equivalent up to combinatorial type.
The collection of all the chambers is the \emph{chamber complex} or
\emph{chamber system} associated with $B$.  Putting all this together
we conclude

\begin{proposition}
  To represent all the possible combinatorial types of
  \emph{non-degenerate} polytopes of the form $P_c = \{x : Bx=c, x
  \geq 0 \}$, for a fixed $B$ and varying vector $c$, it is enough to
  choose one $c$ from each chamber of the chamber complex of $B$.
\end{proposition}

\begin{example}
\label{exm:triprism}
Consider the matrix $B_{2,3}$, the constraint matrix of all $2 \times
3$ transportation polytopes. That is:
\[
B_{2,3}=\left [
\begin {array}{cccccc} 
1&1&1&0&0&0\\
\noalign{\medskip}
0&0&0&1&1&1\\
\noalign{\medskip}
1&0&0&1&0&0\\
\noalign{\medskip}
0&1&0&0&1&0\\
\noalign{\medskip}
0&0&1&0&0&1
\end {array}\right ]
\]
This means, the system $\{B_{2,3}\,y=c,\ y \geq 0\}$ defines the
$2\times 3$ transportation polytopes with marginals $c$.

The columns of the matrix $B_{2,3}$ span a four-dimensional cone in
$\R^5$.  It will be relevant later that if we slice this cone by an
affine hyperplane (such us $\sum y_i =1$) we obtain the
three-dimensional triangular prism shown in Figure~\ref{ejemplo}, but
embedded in $\reals^4$.
 
The chamber complex can be obtained by slicing the prism with the six
planes containing a vertex of the prism and the edge ``opposite'' to
it.  The resulting chamber complex is hard to visualize or draw, even
in this small case, but we will see below how to recover the structure
of the chamber complex for this example using \emph{Gale transforms}.
In particular, as we will see, this decomposes the triangular prism
into 18 chambers.
\end{example}

It is very easy to ``sample'' inside the chamber complex and find
chambers of different numbers of bases, i.e., transportation polytopes
with different number of vertices. One can simply throw random
positive values to the cell entries of an $l \times m \times n$ table
and then compute the 1-marginals or 2-marginals associated to it. But
with this method it is not obvious how to guarantee that one has
obtained all the possible chambers. For this we use the approach based
on Gale transforms and regular triangulations, that we now explain.

A vector configuration $A$ of $r$ vectors in $\R^{(r - d)}$ is called
a \emph{Gale transform} of another vector configuration $B$ of $r$
vectors in $\R^{d}$ if the row space of the matrix with columns given
by $A$ is the orthogonal complement in $\R^r$ of the row space of the
matrix with columns given by $B$. Gale transforms are essential tools
in the study of convex polytopes because the combinatorial properties
of $B$ and $A$ are intimately related (see Chapter 6 in \cite{ziegler}
for details). Proofs of the following statements can be found in
\cite{billeraetal,GKZ}. See also Chapters 4 and 5 in \cite{DRS}.

\begin{lemma} \label{secondary}
\begin{itemize}
\item The chambers of $B$ are in bijection with the regular
  triangulations of the Gale transform $\hat{B}$ of $B$: From a
  chamber in $B$, one can recover a regular triangulation of $\hat{B}$
  via complementation, namely for a basis $\sigma$ of vectors in $B$
  the elements of $\hat{B}$ not belonging to $\sigma$ form a basis for
  $\hat{B}$.  The collection of those bases gives a triangulation of
  $\hat{B}$ (see below for an example).
\item There exists a polyhedron, the \emph{secondary polyhedron},
  whose vertices are in bijection with the regular triangulations of
  the Gale transform $\hat{B}$.
\item The face lattice of the chamber complex of the vector
  configuration $B$ is anti-isomorphic to the face lattice of the
  secondary polyhedron of the Gale transform $\hat{B}$ of $B$. The
  latter is, in turn, isomorphic to the refinement poset of all
  regular subdivisions of $\hat{B}$.
\item If $B$ defines a pointed polyhedral cone (for example, if all
  its entries are non-negative as it is the case for transportation
  polytopes), then its Gale transform $\hat{B}$ is a totally cyclic
  vector configuration. That is, the cone spanned by $\hat{B}$ is the
  whole of $\R^n$.
\end{itemize}
\end{lemma}

\begin{example}[Example~\ref{exm:triprism} continued]
  We take again the matrix $B_{2,3}$ (the case of $2 \times 3$
  transportation polytopes). A Gale transform consists of the columns
  of the matrix
\[
\hat{B}_{2,3}=\left [\begin {array}{cccccc} 1&-1&0&-1&1&0\\\noalign{\medskip}1&0&-1&
-1&0&1\end {array}\right ].
\]

In Figure~\ref{ejemplo} we represent the Gale diagram $\hat{B}_{2,3}$
and its $18$ regular triangulations, each one providing a
combinatorial type of non-degenerate $2\times 3$ transportation
polytope, although repeated combinatorial types occur. The chamber
adjacency, which corresponds to bistellar flips, is indicated by
dotted edges.

\begin{figure}[hbt]
  \begin{center}
    \includegraphics[scale=.25]{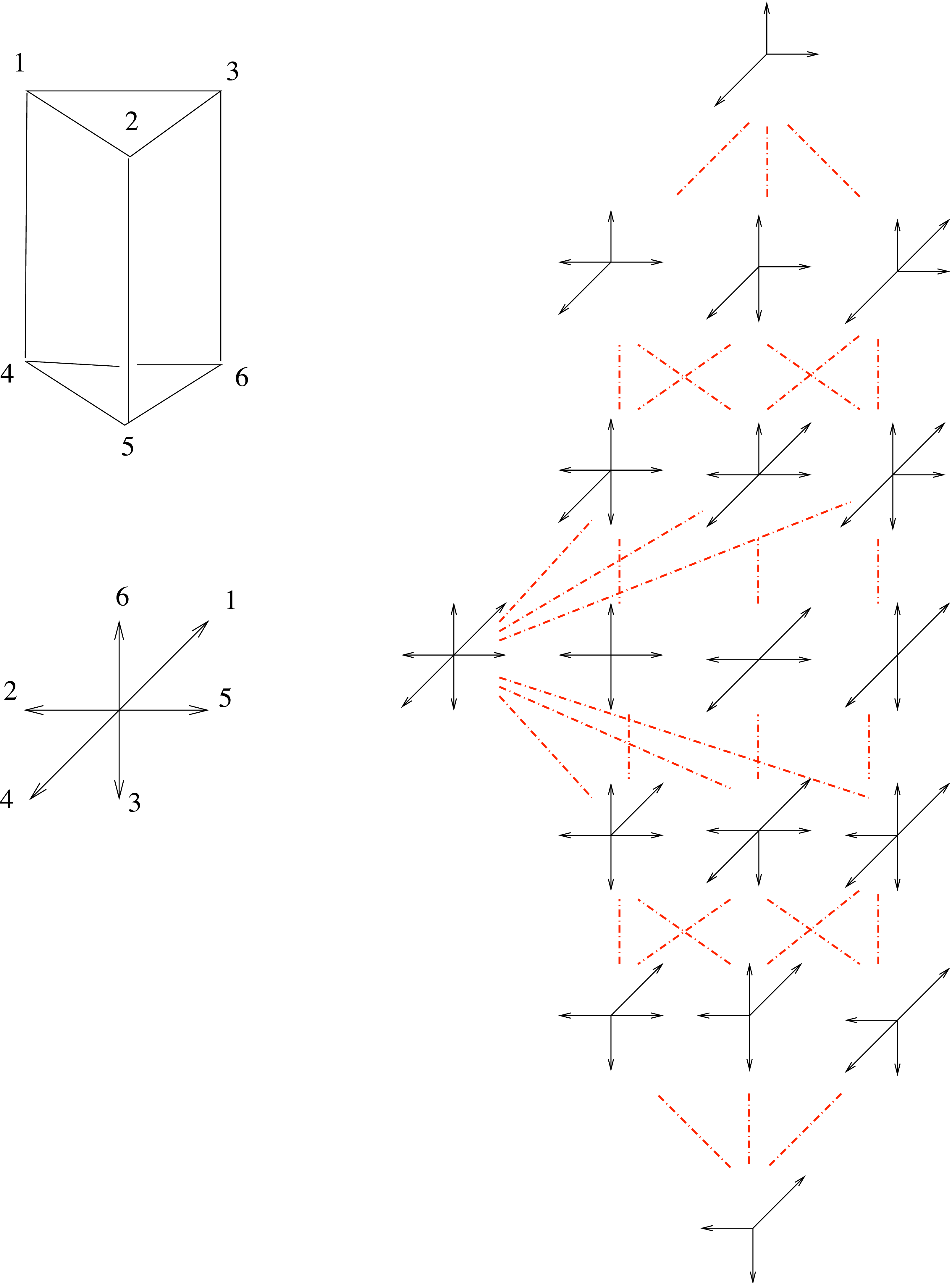}
    \caption{ The 18 regular triangulations of $\hat{B}_{2,3}$, which
      correspond to the 18 combinatorial types of non-degenerate $2 \times 3$
      transportation polytopes.}
       \label{ejemplo}
  \end{center}
\end{figure}
\end{example}

Thus, generating all the combinatorial types of non-degenerate
transportation polytopes is the same as listing the distinct regular
triangulations of the Gale transform of the defining matrix $B$. Note
that in our case $B$ depends only on the sizes $l,m,n$ and the type,
axial or planar, of transportation polytopes we look at. That is, we
have one $B$ for each row of
Tables~\ref{classical_vertices},~\ref{planar_vertices}
and~\ref{axial_vertices}.

Now, it is well-known that the regular triangulations of a vector
configuration can all be generated by applying bistellar flips to a
seed regular triangulation (see
\cite{billeraetal,DRS,ziegler}). Bistellar flips are combinatorial
operations that transform one triangulation into another and
regularity of triangulations can be determined by checking feasibility
of a certain linear program (in our case, the very one that defines
$P_c$). An example of the linear programming feasibility problem is
given in Chapter 5 of \cite{DRS}.

\begin{example}[Example~\ref{exm:triprism} continued]
\label{exm:triprism2b}
Consider the only triangulation of $\hat{B}_{2,3}$ with six cones (the first one of 
the middle row in Figure~\ref{ejemplo}). The necessary and sufficient
conditions in the  non-negative vector $(c_1,c_2,\dots,c_6$) in order to 
produce this triangulation are that each $c_i$ be smaller than the sum
of the two adjacent to it. That is,
\begin{eqnarray*}
c_1 < c_5 + c_6, & c_2 < c_4 + c_6, & c_3 < c_4 + c_5, \cr
c_4 < c_2 + c_3, & c_5 < c_1 + c_3, & c_6 < c_1 + c_2. 
\end{eqnarray*}
Thus, these conditions on the marginals characterize 
the $2\times 3$ transportation polytopes that are hexagons. 
\end{example}

To implement this method we have written a $C^{++}$ program that is
available from the web page of the second author. This program calls
{\tt TOPCOM} (see \cite{topcom}), a package for triangulations that
computes, among other things, the list of all regular triangulations
of a configuration.  Our program also calls {\tt polymake} (see
\cite{polymake}) for the Gale transform, and computes one vector $c$
per chamber. The output is a list of transportation polytopes, one per
chamber, given in the {\tt polymake} file format.

\subsection{Lower and upper bounds via integer programming}

Even for seemingly small cases, such as $3 \times 3 \times 3$
transportation polytopes, listing all chambers (and thus all
combinatorial types of transportation polytopes) is practically
impossible. In these cases we have followed a different approach to at
least obtain upper and lower bounds for the number of vertices of
transportation polytopes. By the discussion above, this is the same as
finding bounds for the number of simplices in triangulations of the
Gale transform $\hat{B}$.  Here we follow the method proposed in
\cite{dhss}, based on the \emph{universal polytope}. This universal
polytope, introduced by Billera, Filliman and Sturmfels in~\cite{bfs},
has \emph{all} triangulations (regular or not) of a given vector
configuration $A$ in $\R^n$ as vertices, and projects to the secondary
polytope. The universal polytope has much higher dimension than the
secondary polytope, in fact its ambience dimension is the number of
possible bases of the configuration $A$, thus no more than $|A|
\choose n+1$.  It has the advantage that the number of simplices in
different triangulations is given by the values of a linear functional
$\psi$.

More precisely, we think of the chambers of $cone(B)$ as the vertices
of the following high-dimensional $0/1$-polytope: Assume $B$ is a
vector configuration with $n$ vectors inside $\R^d$.  Let $N$ be the
number of $d$-dimensional simple cones in $B$.  We define $U_B$ as the
convex hull in $\R^N$ of the set of incidence $0/1$ vectors of all
chambers of $B$. For a chamber $T$ the {\em incidence vector} $v_T$
has coordinates $\,(v_T)_\sigma = 1$ if the basis $\sigma \in T$ and
$\,(v_T)_\sigma = 0$ if $\sigma$ is not a basis of $T$.  The polytope
$U_B$ is the {\em universal polytope} defined in general by Billera,
Filliman and Sturmfels in \cite{bfs} (although there it is defined in
terms of the triangulations of the Gale transform of $B$).

In \cite{dhss}, it was shown that the vertices of the universal
polytope of $B$ are exactly the integral points inside a polyhedron
that has a simple inequality description in terms of the oriented
matroid of $B$ (see \cite{dhss,ziegler} for information on oriented
matroids). The concrete integer programming problems in question were
solved using {\em C-plex Linear Solver}$^{TM}$.  The program to
generate the linear constraints is a small $C^{++}$ available from the
web page of the first author (see \cite{UG}).

\begin{example}[Example~\ref{exm:triprism} continued]
\label{exm:triprism3}

Continuing with the running example, if $B$ is $B_{2,3}$ from above,
then $U_B$ is defined in $\R^{15}$, where each coordinate is indexed
by a $2$-subset $\sigma$ of $\{1,\ldots,6\}$.  Let $S$ denote the set
of all bases.  Then $N = |S| = 15$.  Thus $U_B$ is the convex hull in
$\R^N$ of the incidence vectors $v_T$ corresponding to the $18$
chambers of $B$.  By Lemma \ref{secondary}, this is equivalent to the
convex hull of the incidence vectors $v_T$ of the $18$ triangulations
of $\hat{B}$.  For example, the triangulation $T = \{ \{1,2\},
\{1,3\}, \{2,4\}, \{3,4\} \}$ in Figure \ref{ejemplo} gives the
incidence vector $v_T = e_{\{1,2\}} + e_{\{1,3\}} + e_{\{2,4\}} +
e_{\{3,4\}}$ (where $e_\sigma$ is the basis unit vector in the
direction $\sigma$) as one of the vectors of the convex hull.  

The convex hull of these $18$ incidence vectors is a $6$-dimensional $0/1$
polytope $U_B$ in $\R^{15}$. That the dimension is (at most) six follows from the following considerations:
\begin{itemize}
\item
Since the pairs $\{1,4\}$, $\{2,5\}$ and
$\{3,6\}$ are not full-dimensional and thus never appear as a simplex
in any triangulation $T$ of $\hat{B}$, $U_B$ is contained in the
subspace $x_{\{1,4\}} = x_{\{2,5\}} = x_{\{3,6\}} = 0$.  
\item
Since the vector $1$ has $2$ and $6$ on one side and $5$ and $3$ on the other,
in every triangulation the sum  $x_{\{1,2\}} + x_{\{1,6\}}$ equals the sum 
$x_{\{1,3\}} + x_{\{1,5\}}$ (and it equals zero or one depending on whether
the triangulation uses the vector $1$ or not). This implies the first of the 
following equalities, the rest being the analogue statement for the other five vectors.
\begin{eqnarray*}
x_{\{1,2\}} + x_{\{1,6\}} - x_{\{1,3\}} - x_{\{1,5\}} = &0\\
x_{\{2,3\}} + x_{\{2,4\}} - x_{\{1,2\}} - x_{\{2,6\}} = &0\\
x_{\{1,3\}} + x_{\{3,5\}} - x_{\{2,3\}} - x_{\{3,4\}} = &0\\
x_{\{3,4\}} + x_{\{4,5\}} - x_{\{2,4\}} - x_{\{4,6\}} = &0\\
x_{\{1,5\}} + x_{\{5,6\}} - x_{\{3,5\}} - x_{\{4,5\}} = &0\\
x_{\{2,6\}} + x_{\{4,6\}} - x_{\{1,6\}} - x_{\{5,6\}} = &0
\end{eqnarray*}
Observe that one of these equations is redundant, since the sum 
of the left-hand sides is already zero.

\item Since every triangulation needs to cover the angle between, 
for example, vectors $1$ and $6$, and this angle is covered only by the cones $16$, $12$ and $56$
(see Figure~\ref{ejemplo}), we have that
\[x_{\{1,6\}} + x_{\{1,2\}} + x_{\{5,6\}} =  1.\]
\end{itemize}
The results in~\cite{dhss} say that $U_B$ is the convex hull of the 
non-negative integer points in $\R^{15}$ satisfying this list of equations.

We now denote by $\psi \in \left(\R^N\right)^*$ the cost vector defining
the linear function
\[\psi(x) = (1,1,\ldots,1) \cdot x = \sum_{\sigma \in S} x_\sigma\]
Then the values of $\psi(x)$ on $U_B \cap \{0,1\}^N$ are the only
possible values for the number $f_0$ of vertices of non-degenerate
polytopes of the form $P_c = \{x | Bx=c, x \geq 0 \}$.  In particular,
the solutions to the linear programming relaxations: ``minimize
(respectively maximize) $\psi(x)$ subject to $x \in U_B$'' give lower
(respectively upper) bounds to the possible values for the number
$f_0$ of vertices of non-degenerate polytopes of the form $P_c = \{x |
Bx=c, x \geq 0 \}$.  In the running example, $3 \leq \psi(x) \leq 6$
whenever $x \in U_B \cap \{0,1\}^N$.  From Table
\ref{classical_vertices}, we observe that the number $f_0$ of vertices
of a non-degenerate $2 \times 3$ transportation polytope equals $3$, $4$,
$5$ or $6$.
\end{example}

\begin{example} \label{exm:notBirkhoff} 
Here is an application of our method.
Table~\ref{table:notBirkhoff} is an explicit vector of $2$-marginals
for a $3 \times 3 \times 3$ transportation polytope which has more
vertices (270 vertices) than the generalized Birkhoff polytope, with
only 66 vertices.

\begin{table}[h]
\begin{tabular}{|ccc| |ccc| |ccc|} \hline
      164424  &    324745  &    127239 &   163445   &    49395 &     403568  & 184032  &    123585 &     269245 \\ \hline
      262784  &    601074  &   9369116 &   1151824   &   767866 &    8313284 & 886393  &   6722333 &     935582 \\ \hline
      149654  &   7618489  &   1736281 &   1609500   &  6331023 &    1563901 & 1854344 &    302366 &    9075926 \\ \hline
\end{tabular}
\caption{Counterexample to \cite[open problem 37]{YKK}\label{table:notBirkhoff} }
\end{table}
\end{example}

Based on the data collected from the enumeration process, we also
conjecture to be true:

\begin{conjecture}
  The graph of every non-degenerate $m \times n$ transportation
  polytope has a Hamiltonian cycle ($mn > 4$).
\end{conjecture}

\begin{conjecture}
  If $P$ is a non-degenerate $l \times m \times n$ axial
  transportation polytope ($l, m, n \geq 3$), then the diameter of its
  graph $G(P)$ is equal to $f - d$, where $d = lmn - l - m - n + 2$ is
  the dimension and $f$ is the number of facets of $P$.
\end{conjecture}

For $m \times n$ classical transportation polytopes ($m,n \leq 5$),
there are non-degenerate polytopes where the diameter of the graph
$G(P)$ is strictly less than $f - d$.

\section{Proofs of Theorems~\ref{thm:gcd} and~\ref{thm:22n}}
\label{sec:gcd}

We start with Theorem~\ref{thm:22n}: \emph{The $2 \times 2 \times n$
  planar transportation polytopes are linearly isomorphic to the $2
  \times n$ classical transportation polytopes.}

\begin{lemma}\label{plan22n}
The planar $2 \times m \times n$ transportation polytopes are exactly the  
$m \times n$ transportation polytopes with bounded entries.
\end{lemma}

\bpr 
Every planar $2 \times m \times n$ transportation polytope
\[P\ =\ \left\{(a_{i,j,k})\in\R_+^{2\times m\times n}:
\sum_k a_{i,j,k}=x_{i,j},\ \sum_j a_{i,j,k}=y_{i,k},\ a_{1,j,k}+a_{2,j,k}=z_{j,k}\right\}\]
is linearly isomorphic to an $m \times n$ transportation polytope with bounded entries,
\[Q\ =\ \left\{(a_{1,j,k})\in\R_+^{m\times n}:
\sum_k a_{1,j,k}=x_{1,j},\ \sum_j a_{1,j,k}=y_{1,k},\ a_{1,j,k}\leq z_{j,k}\right\}\ ,\]
via the projection $\R^{2\times m\times n} \rightarrow \R^{m\times n}$ taking $(a_{i,j,k}) \mapsto (a_{1,j,k})$, which maps $P$ bijectively onto $Q$.
Conversely, every $m\times n$ transportation polytope $Q$ with bounded entries
is linearly isomorphic to a planar $2\times m\times n$ transportation polytope $P$
by defining 
$x_{2,j}:=(\sum_k z_{j,k})-x_{1,j},\ j=1,\dots,m$ and $y_{2,k}:=(\sum_j z_{j,k})-y_{1,k},\ k=1,\dots,n$.
\epr

\bprof{Theorem~\ref{thm:22n}} 
Consider any planar $2\times 2\times n$ transportation polytope
\[
P\ =\ \left\{(a_{i,j,k}) \in\R_+^{2\times 2\times n} : \sum_k a_{i,j,k}=x_{i,j},\
a_{i,1,k}+a_{i,2,k}=y_{i,k},\ a_{1,j,k}+a_{2,j,k}=z_{j,k} \right\}\ .
\]
The equations of the last two types imply that for each $k$ we can express all the $a_{i,j,k}$ in
terms of $a_{1,1,k}$ as follows:
\begin{eqnarray*}
a_{1,2,k} &=& y_{1,k} - a_{1,1,k}, \cr
a_{2,1,k} &=& z_{1,k} - a_{1,1,k}, \\
a_{2,2,k} &=& a_{1,1,k} + z_{2,k} - y_{1,k} = a_{1,1,k} + y_{2,k} - z_{1,k}.
\end{eqnarray*}
In particular, $P$ is linearly isomorphic to its projection
\[
Q \ =\ \left\{(a_{1,1,k})  \in\R_+^{2\times 2\times n}:  \alpha_k \le a_{1,1,k}\le \beta_k, \sum_k a_{1,1,k}=x_{1,1} \right\}\ ,
\]
where $\alpha_k= \max\{0, z_{1,k} - y_{2,k}\}=\max\{0, y_{1,k} -
z_{2,k}\}$ and $\beta_k=\min \{y_{1,k},z_{1,k}\}$. Now, by applying a
translation to $Q$, there is no loss of generality in assuming that
$\alpha_k=0$ for all $k$. Then $Q$ is a 1-way transportation polytope
with bounded entries, isomorphic (by Lemma~\ref{plan22n}) to a
$2\times n$ transportation polytope.

Conversely, any $2\times n$ transportation polytope
\[
Q\ =\ \left\{(a_{j,k})\in\R_+^{2\times n}: \sum_k a_{j,k}=x_{j},\
a_{1,k} + a_{2,k}=y_{k}\right\}\ ,
\]
is linearly isomorphic to the following planar $2\times 2\times n$
transportation polytope:
\[
P\ =\ \left\{(a_{i,j,k}) \in\R_+^{2\times 2\times n}: 
\begin{tabular}{l}
$\sum_k a_{1,j,k}=\sum_k a_{2,3-j,k}=x_{j}$, \cr
$\sum_i a_{i,1,k}= \sum_i a_{i,2,k}= 
\sum_j a_{1,j,k}= \sum_j a_{2,j,k}=y_{k}$
\end{tabular}
\right\}\ .
\]
The equations relating the solutions of $Q$ to those of $P$ are
$a_{j,k}=a_{1,j,k}=a_{2,3-j,k}$.\epr One final comment. The above
result is best possible since the list of $2 \times 3 \times 3$ planar
transportation polytopes presented in Table \ref{planar_vertices} is
not the same as the list of $3 \times 3$ classical transportation
problems presented in Table \ref{classical_vertices}.

We now move to Theorem~\ref{thm:gcd}: \emph{The number of vertices of
  a non-degenerate $m \times n$ classical transportation polytope is
  divisible by $\hbox{GCD}(m,n)$.} The first observation, already
hinted in Example~\ref{exm:triprism}, is that the vector configuration
$B_{m,n}$ associated to these transportation polytopes is (a cone
over) the set of vertices of the product $\Delta_{m,n}$ of two
simplices of dimensions $m-1$ and $n-1$. So, we are interested in the
cardinalities of chambers in the product of two simplices. Here and in
what follows we call the \emph{cardinality} of a chamber $c$ of
$B_{m,n}$ the number of bases of $B_{m,n}$ that contain the chamber
$c$. We denote it by $|c|$. The proof of Theorem~\ref{thm:gcd}
consists of the following two steps, which are established
respectively in the two lemmas below:

\begin{itemize}
\item There is a ``seed'' chamber in $\Delta_{m,n}$ whose cardinality
  is indeed a multiple of $\hbox{GCD}(m,n)$.

\item The difference in the cardinalities of any two adjacent chambers
  of $\Delta_m\times \Delta_n$ is a multiple of $\hbox{GCD}(m,n)$.
\end{itemize}

Since the chamber complex is a connected polyhedral complex (where two
adjacent chambers are divided by a hyperplane supported on the vector
configuration) the two lemmas settle the proof.

Let us define the \emph{lexicographic chamber} of $\Delta_{m,n}$
recursively as the (unique) chamber incident to the lexicographic
chamber of $\Delta_{m,n-1}$. The recursion starts with $\Delta_{m,1}$,
which is an $(m-1)$-simplex and contains a unique chamber.  Observe
that the definition of the lexicographic chamber is not symmetric in
$m$ and $n$.  For example, the lexicographic chamber of the triangular
prism $\Delta_{3,2}$ is the one incident to a basis of the prism, and
has cardinality $3$.  The lexicographic chamber of $\Delta_{2,3}$ is
incident to one of the edges parallel to the axis of the prism, and
has cardinality four.

\begin{lemma}
\label{lemma:lex-chamber}
The lexicographic chamber of $\Delta_{m,n}$ is contained in exactly
$m^{n-1}$ simplices.
\end{lemma}

\bpr The cardinality of the lexicographic chamber of $\Delta_{m,n}$
equals the cardinality of the lexicographic chamber of
$\Delta_{m,n-1}$ times the number of vertices of $\Delta_{m,n}$ not
lying in its facet $\Delta_{m,n-1}$. The latter equals $m$.  \epr

When moving from a chamber $c_-$ to an adjacent one $c_+$ we ``cross"
a certain hyperplane $\mathcal{H}$ spanned by all except one of the
elements of any basis containing $c_+$ but not containing $c_-$.  Let
us denote by $C_+$ and $C_-$ the subsets of $B_{m,n}$ lying in the
sides of $\mathcal{H}$ containing $c_+$ and $c_-$ respectively
(Remember that, in our case, $B_{m,n}$ equals the set of vertices of
$\Delta_{m,n}$).
\begin{figure}[hbt]
  \begin{center}
    \includegraphics[scale=0.5]{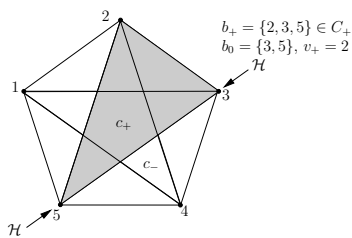}
    \caption{A cross-section of the chamber complex of some cone with two adjacent chambers} \label{cutlemma}
  \end{center}
\end{figure}
Observe also that the common boundary $c_0 \subset
\mathcal{H}$ of $c_+$ and $c_-$ is a chamber in the vector configuration
$B_{m,n}\cap \mathcal H$.

\begin{lemma}
\label{lemma:chambertochamber}
\begin{enumerate}
\item $|c_+|- |c_-| = |c_0| (|C_+|-|C_-|)$.
\item If $B_{m,n}$ is the set of vertices of $\Delta_{m,n}$, 
 then $|C_+|-|C_-|$ is a multiple of $\hbox{GCD}(m,n)$.
\end{enumerate}
\end{lemma}

\bpr A basis $b_+$ contains $c_+$ but not $c_-$ if and only if $b_+$
is of the form $b_0\cup \{v_+\}$, where $b_0$ is a basis of $B\cap
\mathcal H$ containing $c_0$ and $v_+$ is an element of $C_+$.  This
and the analogous property for $c_-$ proves the first part.

For the second part, we restate a few facts in the terminology of
oriented matroids. This makes the proof easier to write (for details
we recommend \cite{redbook}):

\begin{itemize}
\item In oriented matroid terminology a pair $(C_+,C_-)$ consisting of
  the subconfigurations on one and the other side of a hyperplane
  $\mathcal H$ spanned by a subset of $B_{m,n}$ is called a
  \emph{cocircuit} of $B_{m,n}$. That is, part 2 is a statement about
  the cocircuits in the oriented matroid $\mathcal M_{m,n}$ associated
  to the vertices of the product of two simplices.

\item The oriented matroid $\mathcal M_{m,n}$ coincides with the one
  associated to the complete directed bipartite graph $K_{m,n}$.
  (i.e., the complete bipartite graph with all of its edges oriented
  from one part to the other). Thus, part 2 is a statement about the
  cocircuits in the oriented matroid of the directed $K_{m,n}$.

\item The cocircuits of a directed graph $G=(V,E)$ are all read off
  from \emph{cuts} in the graph. By this we mean that the vertex set
  $V$ is decomposed into two parts $(V_+,V_-)$. The cocircuit
  $(C_+,C_-)$ associated to the cut $(V_+,V_-)$ has $C_+$ consisting
  of all the edges directed from $V_+$ to $V_-$ and $C_-$ consisting
  of all the edges directed from $V_-$ to $V_+$.

\end{itemize}

Using the dictionary between the directed graph $K_{m,n}$ and the
product of simplices we can finish the proof.  Let $(V_+, V_-)$ be a
cut in the complete directed bipartite graph $K_{m,n}$. Since our
graph is bipartite, we have $V_+$ and $V_-$ naturally decomposed as
$V_+^{(m)}\cup V_+^{(n)}$ and $V_-^{(m)}\cup V_-^{(n)}$, respectively.
The sizes of $C_+$ and $C_-$ are then:
\[
|C_+| = |V_+^{(m)}| \cdot  |V_-^{(n)}|
\qquad \hbox{and}\qquad
|C_-| = |V_-^{(m)}| \cdot  |V_+^{(n)}|.
\]
Now, using that $|V_+^{(m)}| + |V_-^{(m)}| =m$ and $|V_+^{(n)}| + |V_-^{(n)}| =n$ we get:
\[
|C_+| - |C_-| = 
|V_+^{(m)}| \cdot (n - |V_+^{(n)}|) - |V_+^{(n)}| \cdot (m - |V_+^{(m)}|) = 
|V_+^{(m)}| \cdot n - |V_+^{(n)}| \cdot m,
\] 
which is clearly a multiple of $\hbox{GCD}(m,n)$.
\epr

\section{The diameter of $3$-way axial transportation polytopes}
\label{axialdiam}

Here we consider a 3-way axial transportation polytope $T_{x,y,z}$
defined by certain $1$-marginal vectors $x$, $y$ and $z$.  Recall that
for bounding its diameter there is no loss of generality in assuming
$T_{x,y,z}$ non-degenerate, that is, that $x$, $y$ and $z$ are
sufficiently generic.  In the non-degenerate case, at every vertex $V$
of our polytope exactly $lmn-l-m-n+2$ variables are zero, and exactly
$l+m+n-2$ are non-zero. The set of triplets $(i,j,k)$ indexing
non-zero variables will be called the \emph{support} of the vertex
$V$.

We say that a vertex $V$ of $T_{x,y,z}$ is \emph{well-ordered} if the
triplets $(i,j,k)$ that form its support are totally ordered with
respect to the following \emph{coordinate-wise partial order}:

\begin{equation} \label{order}
(i,j,k)\le (i',j',k') \quad\hbox{ if } \qquad i\le i', \ 
\hbox{and } j\le j', \ \hbox{and}\ k\le k'.
\end{equation}
Observe that a set of $l+m+n-2$ triplets satisfying this must contain
exactly one triplet $(i,j,k)$ with $i+j+k= p$ for each $p=3,\dots,
l+m+n$. Actually, supports of well-ordered vertices are the
\emph{monotone staircases} from $(1,1,1)$ to $(l,m,n)$ in the $l
\times m \times n$ grid.

\begin{lemma}
\label{lemma:well-ordered}
If $x$, $y$ and $z$ are generic, then $T_{x,y,z}$ has a unique
well-ordered vertex $\hat{V}$.
\end{lemma}

\bpr Existence is guaranteed by the ``northwest corner rule
algorithm'', which fills the entries of the table in the prescribed
order (see survey \cite{queyrannespieksma} or exercise 17 in Chapter
six of \cite{YKK}). More explicitly: let $\hat{V}_{l,m,n}=
\min\{x_l,y_m,z_n\}$. Genericity implies that the three values $x_l$,
$y_m$ and $z_n$ are different. Without loss of generality we assume
that the minimum is $z_n$. Then, our choice of $\hat{V}_{l,m,n}$ makes
$\hat{V}_{i,j,n}=0$ for every other pair $(i,j)$. The rest of our
vertex $\hat{V}$ is a vertex of the $l\times m\times (n-1)$ axial
transportation polytope with margins $x'=(x_1, \dots, x_{l-1}, x_l -
z_n)$, $y'=(y_1, \dots, y_{m-1}, y_m - z_n)$, and $z'=(z_1, \dots,
z_{n-1})$.
  
Uniqueness follows from the same argument, simply noticing that the
support of a well-ordered vertex always contains the entry $(l,m,n)$,
and no other entry from one of the three planes $(l,*,*)$, $(*,m,*)$
and $(*,*,n)$. This, recursively, implies that the vertex can be
obtained by the northwest corner rule.  \epr

\begin{remark}
  Another proof of Lemma~\ref{lemma:well-ordered} can be done using
  the formalism of chambers developed in the previous sections: it is
  obvious (and is proved in~\cite{dhss}) that if $c$ denotes a chamber
  of $B$, and $T$ is a triangulation of $\cone(B)$, then there is a
  unique maximal-dimension simplex in $T$ that contains $c$. Thus,
  Lemma~\ref{lemma:well-ordered} follows from the fact that monotone
  staircases in the $l\times m\times n$ grid form a triangulation of
  the vector configuration $B_{l,m,n}$ of axial $l\times m\times n$
  transportation polytopes.  The latter is well-known, once we observe
  that $B_{l,m,n}$ is the vertex set of a product of three simplices.
  The triangulation in question is called the ``staircase
  triangulation'' of it (see Chapter 6 of \cite{DRS}).
\end{remark}

\begin{example}\label{axial333}
  To illustrate Lemma \ref{lemma:well-ordered} consider the
  non-degenerate $3 \times 3 \times 3$ axial transportation polytope
  $T_{x,y,z}$ with:
\[
\sum_{j,k} a_{1,j,k} = 112\qquad 
\sum_{j,k} a_{2,j,k} = 18\qquad 
\sum_{j,k} a_{3,j,k} = 30\qquad 
\]
\[
\sum_{i,k} a_{i,1,k} = 40\qquad 
\sum_{i,k} a_{i,2,k} = 6\qquad 
\sum_{i,k} a_{i,3,k} = 114\qquad 
\]
\[
\sum_{i,j} a_{i,j,1} = 82\qquad 
\sum_{i,j} a_{i,j,2} = 44\qquad 
\sum_{i,j} a_{i,j,3} = 34\qquad 
\]

The unique well-ordered vertex $\hat{V}$ of $T_{x,y,z}$ has the 
non-zero coordinates
$a_{(1,1,1)}=40$,
$a_{(1,2,1)}=6$,
$a_{(1,3,1)}=36$,
$a_{(1,3,2)}=30$,
$a_{(2,3,2)}=14$,
$a_{(2,3,3)}=4$, and
$a_{(3,3,3)}=30$.
Note that the non-zero entries of $\hat{V}$ are totally ordered (they 
are presented above in increasing order) with respect to (\ref{order}).
Figure \ref{monostair} depicts the associated monotone staircase.
\begin{figure}[hbt]
  \begin{center}
    \includegraphics[scale=0.7]{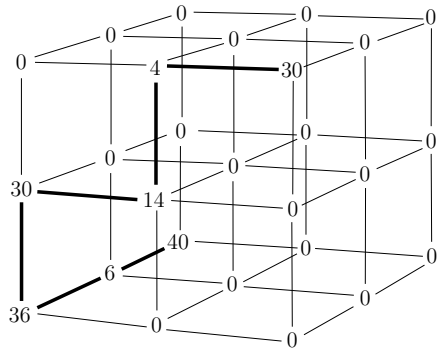}
    \caption{A well-ordered vertex and its staircase.} \label{monostair}
  \end{center}
\end{figure}

\end{example}

Our bound on the diameter of $T_{x,y,z}$ is based on an explicit path
that goes from any initial vertex $V$ of $T_{x,y,z}$ to the unique
well-ordered vertex $\hat{V}$. To build this path we rely on the
following stratified version of the concept of well-ordered vertex.
We say that a vertex $V$ of $T_{x,y,z}$ is \emph{well-ordered starting
  at level $p$}, where $p$ is an integer between 3 and $l+m+n$ if:

\begin{enumerate}
\item For each $q=p,\dots, l+m+n$, the support of $V$ contains exactly
  one triplet $(i,j,k)$ with $i+j+k = q$.
\item Those triplets are well-ordered.  (The partial
  order given in (\ref{order}) is a total order on these triplets.)
\item All other triplets in the support have entries which are
  index-wise smaller than or equal to those of the unique triplet
  $(i_0,j_0,k_0)$ with $i_0+j_0+k_0=p$.
\end{enumerate}

For example, the only vertex ``well-ordered starting at level 3'' is
the well-ordered vertex $\hat{V}$. Slightly less trivially, it is also
the unique vertex ``well-ordered starting at level 4''.  On the other
extreme, all vertices that contain $(l,m,n)$ as a support triplet are
well-ordered starting at level $l+m+n$.  Observe that from any vertex
of $T_{x,y,z}$ we can move, by a single pivot edge in the sense of the
simplex method, to another vertex containing any prescribed entry
$(i,j,k)$ to be non-zero. In particular, we can move to a vertex that
has $(l,m,n)$ in its support.  So, we can assume from the beginning
that $(l,m,n)$ is in the support of our initial vertex $V$, and will
add one to the count of edges traversed to arrive to $\hat{V}$.

\begin{lemma}\label{lemma:decrease-p}
  If $V$ is a vertex of $T_{x,y,z}$ that is well-ordered starting at
  level $p \in\{5,\dots, l+m+n\}$, then there is a path of at most
  $2(p-4)$ edges of $T_{x,y,z}$ that leads from $V$ to a vertex that
  is well-ordered starting at level $p-1$.
\end{lemma}

\bpr Let $(i_0,j_0,k_0)$ be the unique triplet in the support of $V$
with $i_0+j_0+k_0=p$.  We first observe that there is no loss of
generality in assuming that $p=l+m+n$ (that is,
$(i_0,j_0,k_0)=(l,m,n)$). This is because the vertices of $T_{x,y,z}$
that are well-ordered starting at level $p$ and agree with $V$ in all
the triplets with sum of indices greater than or equal to $p$ are the
vertices of a non-degenerate $i_0 \times j_0 \times k_0$ axial
transportation polytope, obtained as in the proof of
Lemma~\ref{lemma:well-ordered}.

So, from now on we assume that $V$ is well-ordered starting at level
$p=l+m+n$.  Let $S_1$ be the set of support triplets in $V$, other
than $(l,m,n)$, that have first index equal to $l$.  Similarly, let
$S_2$ and $S_3$ be the sets of support triplets that have,
respectively, second and third indices equal to $m$ and $n$.

\begin{figure}[hbt]
  \begin{center}
    \includegraphics[scale=0.7]{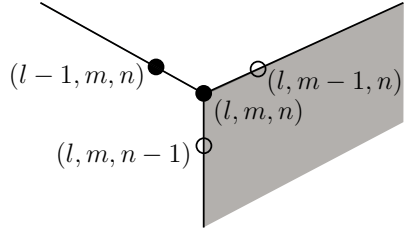}
    \caption{A well-ordered vertex starting at level $p-1$.} \label{lower_p}
  \end{center}
\end{figure}

Our goal is to modify $V$ until one of $S_1$, $S_2$ or $S_3$ becomes
empty, but always keeping the triplet $(l,m,n)$ in the support.  Once
this is done, a single pivot step can be used to obtain a vertex that
is well-ordered starting at level $p-1$ as follows: Without loss of
generality assume that $S_1$ is empty (the cases when $S_2$ or $S_3$
are empty are treated identically). In particular, neither $(l,m-1,n)$
nor $(l,m,n-1)$ are in the support.  If $(l-1,m,n)$ is in the support
then our vertex is already well-ordered starting at level $p-1$. If
not, we do the pivot step that inserts $(l-1,m,n)$. This pivot step
cannot remove $(l,m,n)$ or insert $(l,m-1,n)$ or $(l,m,n-1)$ in the
support.  (The $(l,m,n)$ coordinate is not removed from the support
since the entry remains constant in this pivot.  The $(l,m-1,n)$ and
$(l,m,n-1)$ coordinates remain zero because only non-zero entries of
$V$ and the entry $(l-1,m,n)$ change in the pivot). This pivot
produces a vertex well-ordered starting at level $p-1$.  Figure
\ref{lower_p} gives a picture for this case.

Given a vertex $V$ well-ordered at level $p$, we specify a sequence of
pivots in the graph of $T_{x,y,z}$ to a vertex $V'$ such that one of
$S_1$, $S_2$ or $S_3$ is empty for $V'$.  Lemma~\ref{lemma:decrease-S}
below shows how to get such a $V'$ in a number of steps bounded above by
\[
2|S_1 \cup S_2 \cup S_3| - 3 \le 2 (p-3) - 3 = 2p -9.
\] 
In one more step, that is, at most $2p-8$, we get to a vertex that is
well-ordered starting at level $p-1$. This completes the proof of our
lemma.
\epr

For Lemma~\ref{lemma:decrease-S} let us introduce the following notation:
\[
R_1:= S_1 \setminus (S_2\cup S_3),\quad
R_2:= S_2 \setminus (S_1\cup S_3),\quad
R_3:= S_3 \setminus (S_1\cup S_2),\quad
\]
\[
R_{12}:=S_1\cap S_2, \quad
R_{13}:=S_1\cap S_3, \quad
R_{23}:=S_2\cap S_3.
\]
That is, $R_i$ consists of the elements of $S_1 \cup S_2 \cup S_3$ that belong only to $S_i$,
and $R_{ij}$ of those that belong to $S_i$ and $S_j$. Observe that, by definition, no element belongs
to the three $S_i$'s, so that $S_1 \cup S_2 \cup S_3$ is the disjoint union of these six subsets.

\begin{lemma}
\label{lemma:decrease-S}
With the above notation and the conditions of the proof of the previous lemma, suppose that 
no $S_i$ is empty. Then:
\begin{enumerate}
\item If both $R_i$ and $R_{jk}$ are non-empty, with $\{i,j,k\}=\{1,2,3\}$, then there is a single pivot step that decreases $|S_1 \cup S_2 \cup S_3|$.
\item If the three $R_{ij}$'s are non-empty, then there is a sequence of two pivot steps that decreases $|S_1 \cup S_2 \cup S_3|$.
\item If the three $R_i$'s are non-empty, then there is a sequence of two pivot steps that decreases $|S_1 \cup S_2 \cup S_3|$.
\item If none of the above happens, then $S_1 \cup S_2 \cup S_3$ is contained in one of the $S_i$'s, say $S_1$. Then, there is a sequence of $|S_1|-1$ pivot steps that makes $S_1 \cup S_2 \cup S_3$ empty.
\end{enumerate}
All in all, there is a sequence of at most $2|S_1 \cup S_2 \cup S_3| - 3$ pivot steps that makes some $S_i$ empty.
\end{lemma}

\bpr 
Let us first show how the conclusion is obtained.  We argue by induction on $|S_1 \cup S_2 \cup S_3|$, the base case being $|S_1 \cup S_2 \cup S_3| = 2$, which is the minimum to have no $S_i$ empty. The base case implies we are in the situation of either part (1) or part (4), and a single pivot step makes an $S_i$ empty.

If $|S_1 \cup S_2 \cup S_3| > 2$ and one of the conditions (1), (2) or (3) holds, then we do the step or the two pivot steps mentioned there and apply induction. If none of these three conditions hold then it is easy to see that (4) must hold. (Remember that we are assuming that no $S_i$ is empty, and $S_i=R_i\cup R_{ij}\cup R_{ik}$).
Part (4) guarantees we have a sequence of $|S_1|-1 \le 2|S_1 \cup S_2 \cup S_3| - 3$ pivot steps that makes an $S_i$ empty.

So, let us prove each of the four items in the lemma. Let $\alpha$ denote the entry $(l,m,n)$.

\begin{enumerate} 
\item Suppose without loss of generality that $R_{12}$ and $R_3$ are not empty.
Let  $\beta\in R_{12}$ and $\gamma\in R_3$.
  The reader may find it useful to follow our proof using Figure
  \ref{case1} which depicts the situation.
\begin{figure}[hbt]
  \begin{center}
    \includegraphics[scale=0.7]{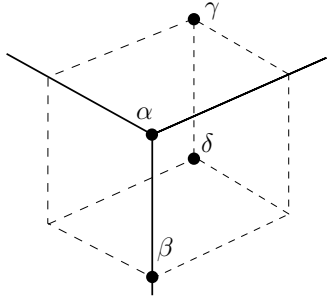}
    \caption{A layout of entries $\alpha$, $\beta$, $\gamma$ and $\delta$.} \label{case1}
  \end{center}
\end{figure}
Observe
that $\alpha=(l,m,n)$ is the index-wise maximum of $\beta$ and
$\gamma$. Let $\delta$ be the index-wise minimum of them. First
observe that $\delta$ is not in the support of vertex $V$. Otherwise,
we could add $\pm
\frac{1}{2}\min\{a_\alpha,a_\beta,a_\gamma,a_\delta,\}
\left(e_\alpha+e_\delta-e_\beta-e_\gamma\right)$ to $V$ and stay in
$T_{x,y,z}$. Hence, $V$ would be a convex combination of two other
points from $T_{x,y,z}$ (and thus not a vertex), parallel to the
direction of $e_\alpha+e_\delta-e_\beta-e_\gamma$ (here $e_{i,j,k}$
denotes the basis unit vector in the direction of the variable
$a_{i,j,k}$).

Next, consider $V' = V + \min\{a_\beta, a_\gamma\}
\left(e_\alpha+e_\delta-e_\beta-e_\gamma\right)$.  Observe that $V'$
has different support than $V$ since either $\beta$ or $\gamma$ has
been removed (not both, because $a_\beta\ne a_\gamma$ by
non-degeneracy).  Also, since $V'$ cannot have support strictly
contained in that of $V$, $\delta$ must have been added. That is, the
supports of the vertices $V$ and $V'$ differ in the deletion and
insertion of a single element, which means they are adjacent in the
graph of the polytope $T_{x,y,z}$. As desired, when going from $V$ to $V'$
the cardinality of $S_1\cup S_2 \cup S_3$ is decreased by one.

\begin{example}[Example \ref{axial333} continued]
Consider the vertex $V$ with non-zero coordinates 
$a_{(1,1,2)}=28$, 
$a_{(2,1,2)}=12$, 
$a_{(2,2,3)}=6$, 
$a_{(1,3,2)}=2$, 
$a_{(1,3,1)}=82$, 
$a_{(3,3,2)}=2$, and 
$a_{(3,3,3)}=28$.
In this example, $\alpha=(3,3,3)$, $\beta=(3,3,2)$, $\gamma=(2,2,3)$, 
and $\delta=(2,2,2)$.  After clearing $a_\beta$, we arrive at the vertex 
$V'$ with non-zero coordinates
$a'_{(1,1,2)}=28$, 
$a'_{(2,1,2)}=12$, 
$a'_{(2,2,3)}=4$, 
$a'_{(1,3,2)}=2$, 
$a'_{(1,3,1)}=82$, 
$a'_{(2,2,2)}=2$, and 
$a'_{(3,3,3)}=28$.
\end{example}

\item Suppose now that none of the $R_{ij}$'s is empty, and let $\beta\in R_{13}$, $\gamma\in R_{23}$ and $\delta' \in R_{12}$.
we apply the same pivot as in case one, which makes $\delta$, the coordinate-wise minimum of $\beta$ and $\gamma$, enter the support. This pivot does not decrease $|S_1\cup S_2 \cup S_3|$, but it leads to a situation where we have $\delta \in R_3$ and 
$\delta' \in R_{12}$. Hence, we can apply part one and decrease $|S_1\cup S_2 \cup S_3|$ with a second step.

\item[4.] Let us prove now part (4) and leave (3), which is more complicated, for the end. Observe that if 
$S_1\cup S_2 \cup S_3 = S_1$ but $S_2$ and $S_3$ are not empty, then necessarily $R_{12}$ and $R_{13}$ are both non-empty. While this holds, we can do the same pivot steps as before with a $\beta\in R_{12}$ and a $\gamma\in R_{13}$. Each step decreases by one the cardinality of $R_{12}\cup R_{13}$, increasing the cardinality of $R_1$. The process finishes when $R_{12}$ (hence $S_2$) or $R_{13}$ (hence $S_3$) becomes empty, which happens, in the worst case, in $|S_1| -1$ steps.

\item[3.] Finally, consider the case where the three $R_i$'s are non-empty.
 Let $\beta = (l,j_1,k_1)\in R_1$,
  $\gamma=(i_2,m,k_2)\in R_2$ and $\delta=(i_3,j_3,n)\in R_3$, and, as
  before, $\alpha=(l,m,n)$.  Figure \ref{case2} depicts the situation
  to help the reader with following our proof.
\begin{figure}[hbt]
  \begin{center}
    \includegraphics[scale=0.7]{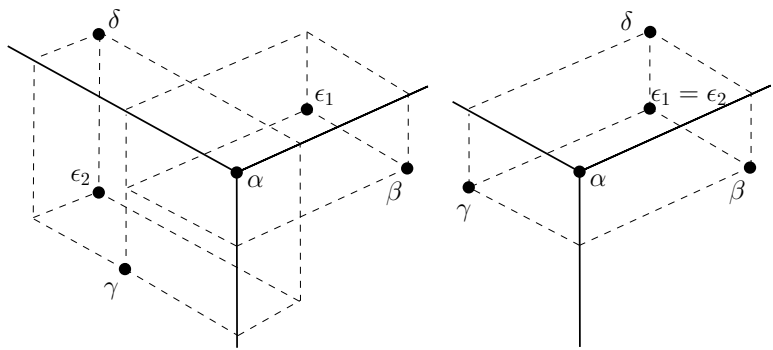}
    \caption{Possible layouts of the entries $\alpha$, $\beta$, $\gamma$, $\delta$, $\epsilon_1$ and $\epsilon_2$.}\label{case2}
  \end{center}
\end{figure}
Let $\epsilon_1$ and $\epsilon_2$ be two triplets of indices with the
property that $\{\alpha,\epsilon_1,\epsilon_2\}$ and $\{\beta, \gamma,
\delta\}$ use exactly the same three first indices, the same three
second indices, and the same three third indices. For example, let us make the
choice $\epsilon_1=(i_2,j_1,k_1)$ and $\epsilon_2=(i_3,j_3,k_2)$ as in
the left part of Figure \ref{case2}.  By non-degeneracy, the smallest
value among $a_\beta$, $a_\gamma$ and $a_\delta$ at $V$ is unique.  We
assume without loss of generality that the smallest among them is
$a_\beta$. Let $W$ be the point of $T_{x,y,z}$ obtained by changing
the following six coordinates:
\[
a'_\alpha = a_\alpha+ a_\beta ,\
a'_\beta = a_\beta - a_\beta=0,\
a'_\gamma = a_\gamma - a_\beta,\
a'_\delta = a_\delta - a_\beta,\
a'_{\epsilon_1} = a_{\epsilon_1}+ a_\beta ,\
a'_{\epsilon_2} = a_{\epsilon_2}+ a_\beta.
\]
It may occur that $\epsilon_1=\epsilon_2$, as shown in the right side
of Figure \ref{case2}. Then we do the same pivot except we increase
the corresponding entry $a_{\epsilon_1}= a_{\epsilon_2}$ twice as
much.

Observe that one of $\epsilon_1$ or $\epsilon_2$ may already be in the
support of $V$, but not both: Otherwise $W$ would have support
strictly contained in that of $V$, which is impossible because $V$ is
a vertex and has minimal support. If one of $\epsilon_1$ or
$\epsilon_2$ were already in the support of $V$, or if
$\epsilon_1=\epsilon_2$, then $W$ is a vertex and we take $V'=W$. As
in the first case, $V'$ is obtained from $V$ by traversing a single
edge and has one less support element in $S_1 \cup S_2 \cup S_3$ than $V$: None of $\epsilon_1$ or $\epsilon_2$
can have a common entry with $\alpha$, since none of $\beta$, $\gamma$
and $\delta$ has two common entries with $\alpha$.

However, if $\epsilon_1$ and $\epsilon_2$ are different and none of
them was in $V$, then $W$ has one-too-many elements in its support to
be a vertex, which means it is in the relative interior of an edge $E$
and $L=VW$ is not an edge.  See Figure \ref{octa}. Moreover, both $E$ and
$VW$ lie in a
two-dimensional face $F$.  
This is so because every support containing the support of a vertex defines
a face of dimension equal the excess of elements it has. In our case, $F$ is 
the face 
with support $\hbox{support}(V) \cup \hbox{support}(W) =
\hbox{support}(V) \cup \{\epsilon_1,\epsilon_2\}$.

\begin{figure}[hbt]
  \begin{center}
    \includegraphics[scale=0.7]{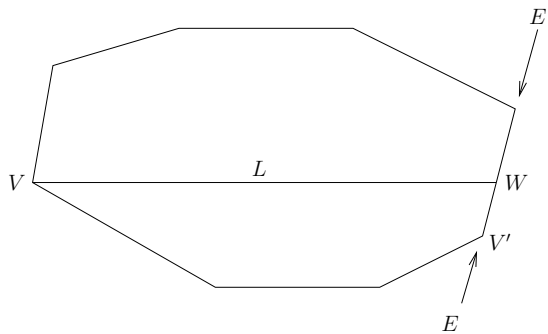}
    \caption{The octagon containing segment $VW$ arising when entries 
$\epsilon_1$ and $\epsilon_2$ are different and none of
them was in the support of $V$} \label{octa}
  \end{center}
\end{figure}

We now look more closely at the structure of $F$. Each edge $H$ of $F$
is the intersection of $F$ with a facet of our transportation polytope. That is,
there is a unique variable $\eta$ that is constantly zero along $H$  but not zero as we
move on $F$ in other directions. For example, since $\epsilon_1$ and $\epsilon_2$
are zero at $V$ but not constant on $F$ (they increase along $L$),  $V$
is the common end of the edges defined by  $\epsilon_1$ and $\epsilon_2$.

Our goal is to show that there is a vertex $V'$ of $F$ at distance at most two from $V$
and incident to the edge defined by one of the variables  $\beta$, $\gamma$ and $\delta$.
At such a vertex $V'$ we will have decreased $|S_1\cup S_2\cup S_3|$ by one, as claimed.
The key remark is that there are at most two edges of $F$ not produced by one of the variables
$\alpha$, $\beta$, $\gamma$, $\delta$, $\epsilon_1$, and $\epsilon_2$:
every variable $\eta$ other than those six is constant along $L$, so it either produces an edge parallel to $L$ or no edge at all. In particular, $F$ is at most an octagon, as in Figure~\ref{octa}. Now:

\begin{itemize}
\item If $F$ has  five or less edges, then every vertex of $F$ is at distance one or two from $V$.
Take as $V'$ either end of the end-point $W$ of $L$. This works because at $W$ one of $\beta$, $\gamma$
or $\delta$ is zero, by construction.

\item If $F$ has six or more edges, then the two vertices $V'$ and $V''$ of $F$ at distance two from $V$ are at distance at least two from each other. So, together they are incident to four different edges, none of which is the edge of $\epsilon_1$ or $\epsilon_2$. (Remember that the edges of $\epsilon_1$ and $\epsilon_2$ are incident to $V$).
At least one of these four edges is defined by $\beta$, $\gamma$, or $\delta$,
because there are (at most) three other possible edges: the one of $\alpha$ and two parallel to $L$.
\end{itemize}
\epr

\end{enumerate}

\begin{example}
To make ideas completely clear, using the same polytope $T_{x,y,z}$ as
in Example \ref{axial333}, we consider its vertex $V$ with non-zero coordinates
$v_{(1,1,3)}=25$, 
$v_{(3,1,1)}=15$, 
$v_{(3,2,1)}=6$, 
$v_{(1,3,1)}=61$, 
$v_{(1,3,2)}=26$, 
$v_{(2,3,2)}=18$ and 
$v_{(3,3,3)}=9$.
Here, $\alpha=(3,3,3)$, $\beta=(3,2,1)$, $\gamma=(2,3,2)$, $\delta=(1,1,3)$, 
$\epsilon_1=(2,2,1)$, and $\epsilon_2=(1,1,2)$.
The triplet $\beta$ is not in the support of $W$, and $W$
has non-zero coordinates
$w_{(1,1,3)}=19$, 
$w_{(1,1,2)}=6$,
$w_{(3,1,1)}=15$, 
$w_{(2,2,1)}=6$, 
$w_{(1,3,1)}=61$, 
$w_{(1,3,2)}=26$, 
$w_{(2,3,2)}=12$, and 
$w_{(3,3,3)}=15$.

The vertices of $T_{x,y,z}$ with support contained in 
$support(V) \cup \{\epsilon_1, \epsilon_2\}$ form the $4$-gon $F=VBCD$ where
$B$ is the vertex with non-zero coordinates
\[
b_{(1,1,3)}=12, 
b_{(1,1,2)}=26, 
b_{(3,1,1)}=2, 
b_{(3,2,1)}=6, 
b_{(3,3,3)}=22, 
b_{(2,3,2)}=18,
b_{(1,3,1)}=74,
\]
$C$ is the vertex with non-zero coordinates
\[
c_{(1,1,3)}=22, 
c_{(3,1,1)}=18, 
c_{(2,2,1)}=6, 
c_{(3,3,3)}=12,
c_{(1,3,2)}=32, 
c_{(2,3,2)}=12, 
c_{(1,3,1)}=58
\]
and $D$ is the vertex with non-zero coordinates
\[
d_{(1,1,3)}=6, 
d_{(1,1,2)}=32, 
d_{(3,1,1)}=2, 
d_{(2,2,1)}=6, 
d_{(3,3,3)}=28, 
d_{(2,3,2)}=12,
d_{(1,3,1)}=74.
\]
Note $W$ is in the edge $E=CD$.  We let $V' = D$, the endpoint of $E$ 
closer to $V$.  Thus, we use one edge to go from $V$ to $V'$.
\end{example}

\bprof{Theorem \ref{theorem:main}} Starting with any vertex
``well-ordered starting at level $p=l+m+n$'' (which can be reached in a
single step) we use Lemma \ref{lemma:decrease-p} to decrease one unit
by one the level at which our vertex starts to be well-ordered until
we reach the unique well-ordered vertex $\hat{V}$. Thus, the number of
steps needed to go from an arbitrary vertex $V$ to $\hat{V}$ is at most
\[
1 + \sum_{q=5}^p 2(q-4) =  1 + 2 \sum_{q=1}^{p-4} q
=1 + 2{p-3 \choose 2} \le (p-3)^2.
\]
To go from one arbitrary vertex to another, twice as many steps suffice.
\epr

\noindent{\bf Remark:} 
The whole proof can be generalized to arbitrary axial $d$-way tables,
instead of $d=3$, without much effort. Everything in
Lemma~\ref{lemma:well-ordered} goes through without change, as well as
the definition of ``well-ordered starting at level $p$''.  In the
other arguments, the first change is that we have $d$ sets $S_1,\dots
, S_d$ instead of just three. In particular, in the proof of
Lemma~\ref{lemma:decrease-S}, the worst case will be that of $d-1$
different $\epsilon$'s, which gives a face of dimension $d-1$. Hence,
the bound given in the statement of Lemma~\ref{lemma:decrease-p} can
be substituted to the maximum diameter of a simple polytope of
dimension $d-1$ with at most $p$ facets.  This still yields a
polynomial bound for any fixed value of $d$. We leave the details for the
interested reader.

\section{Acknowledgements}

The authors are grateful to Cor Hurkens, Fu Liu, Maurice Queyranne,
Leen Stougie and G\"unter Ziegler for useful conversations and
references. The authors are also grateful to the referees of the paper
for their thoughtful comments and remarks.


\end{document}